\documentclass[10pt,A4paper,twoside]{article}
\usepackage{amstext}
\usepackage{theorem}
\usepackage{amssymb}
\usepackage{latexsym}
\usepackage{amsmath}
\usepackage{amscd}

\newtheorem{theorem}{Theorem}[section]
\newtheorem{proposition}[theorem]{Proposition}
\newtheorem{lemma}[theorem]{Lemma}

\newtheorem{conjecture}[theorem]{Conjecture}
\newtheorem{corollary}[theorem]{Corollary}

\newtheorem{definition}[theorem]{Definition}

\newtheorem{problem}[theorem]{Problem}

\newcommand{\Z}{\mathbb{Z}}
\newcommand{\Q}{\mathbb{Q}}
\newcommand{\R}{\mathbb{R}}
\newcommand{\C}{\mathbb{C}}

\newcommand{\Pro}{\mathbb{P}}

\newcommand{\Acal}{\mathcal{A}}

\newcommand{\Kcal}{\mathcal{K}}
\newcommand{\Lcal}{\mathcal{L}}
\newcommand{\Mcal}{\mathcal{M}}

\newcommand{\Ocal}{\mathcal{O}}
\newcommand{\ord}{\mbox{\upshape{ord}}}
\newcommand{\BP}{\mathbf{BP}}
\newcommand{\im}{\mbox{\upshape{im}}}

\newenvironment{proof}{{\it Proof\,:}}{\hfill$\Diamond$}

\title{Representation of squares by monic second degree polynomials in the field of $p$-adic meromorphic functions}

\author{Hector Pasten\\
University of Concepci\'on}
\date{}

\begin{document}

\date{}
\maketitle
\begin{abstract}
We prove a result on the representation of squares by second degree polynomials in the field of $p$-adic meromorphic functions in order to solve positively B\"uchi's $n$ squares problem in this field (that is, the problem of the existence of a constant $M$ such that any sequence $(x_n^2)$ of $M$ - not all constant - squares whose second difference is the constant sequence $(2)$ satisfies $x_n^2=(x+n)^2$ for some $x$). We prove (based on works by Vojta) an analogous result for function fields of characteristic zero, and under a Conjecture by Bombieri, an analogous result for number fields. Using an argument by B\"uchi, we show how the obtained results improve some theorems about undecidability for the field of $p$-adic meromorphic functions and the ring of $p$-adic entire functions.  
\end{abstract}

\tableofcontents
\section{Introduction}

In 1970, after the work developed by M. Davis, H. Putnam and J. Robinson, Hilbert's Tenth Problem was answered negatively by Y. Matiyasevich. In logical terms, it was shown that the positive existential theory of $\Z$ in the language of rings $\Lcal_R=\{0,1,+,\cdot\}$ is undecidable, which means that there exists no algorithm to decide whether a system of diophantine equations (or equivalently, a single diophantine equation) has integer solutions or not. 

Soon after the problem was solved, J. R. B\"uchi proved in an unpublished work (see \cite{Lipshitz}) that a positive answer to a certain problem in Number Theory (which we write here $\BP(\Z)$) would allow to show that there exists no algorithm to decide whether a system of diagonal quadratic diophantine equations has integer solutions or not. 

The number-theoretical problem $\BP(\Z)$ is based on the following observation. If we consider the first difference of a sequence of consecutive integer squares (for example $1,4,9,16$), we obtain a sequence of odd integers (in our example $3,5,7$). Hence, the second difference is constant and equal to two. One may ask whether a sequence of squares having second difference equal to two must be a sequence of consecutive squares. The sequence $6^2,23^2,32^2,39^2$ shows that in general such a reciprocal is not true.

\begin{problem}[$\BP(\Z)$]\label{buchi} Does there exist an integer $M$ such that the following happens:\\
If the second difference of a sequence $(x_i^2)_{i=1}^M$ of integer squares is constant and equal to $2$, then there exists an integer $\nu$ such that $x_i^2=(\nu + i)^2$ for $i=1,\ldots,M$ (that is, the squares must be consecutive).
\end{problem}

This problem became known as the \emph{$n$ Squares Problem} or \emph{B\"uchi's Problem}. Numerical evidence suggests that $M=5$ should work, but $\BP(\Z)$ still is an open problem. 

Assuming a positive answer to $\BP(\Z)$, B\"uchi was able to construct an algorithm to do the following: given a diophantine equation, to construct a system of quadratic diagonal equations such that the former has a solution if and only if the latter has. Therefore, using the negative answer given to Hilbert's tenth problem and assuming a positive answer to $\BP(\Z)$, we get the non-existence of an algorithm to decide whether a system of diophantine diagonal quadratic equations has an integer solution. The non-existence of such an algorithm can be shown to be equivalent to the undecidability of the positive existential theory of $\Z$ over the language $\Lcal_2=\{0,1,+,P_2\}$ where $P_2(x)$ is interpreted as `$x$ is a square'. 

In order to get similar consequences in Logic for other rings of interest, and motivated by the arithmetical interest of the problem, several authors have studied variants of $\BP(\Z)$. A natural thing to do is to replace the ring $\Z$ by another commutative ring $A$ with unit. Depending on the ring, we sometimes need to make additional hypothesis in the statement of $\BP(A)$:
\begin{itemize}
\item If $A$ is a ring of functions of characteristic zero in the variable $z$, then we ask for at least one $x_i$ to be non-constant. 
\item If $A$ is a ring of positive characteristic, then we ask $M$ to be at most the characteristic of $A$.
\end{itemize}

The positive existential $\Lcal_2$-theory of a ring is usually much weaker than its positive existential $\Lcal_R$-theory. But when B\"uchi's problem has a positive answer for a ring $A$ then those theories for $A$ are (in general) equivalent. This is what happens for example for $p$-adic analytic functions and for $p$-adic meromorphic functions (see Section \ref{ResLogic}).

In this work, we will solve $\BP(A)$ for some rings of functions (namely, the field of p-adic meromorphic functions and function fields of curves in characteristic zero) by showing in each case a somewhat stronger result on representation of squares by polynomials, in the spirit of the following:

\textit{Given an algebra $A\subseteq B$, there exists a constant $M$ satisfying the following condition:\\
For any set $\{a_1,\ldots,a_M\}$ of $M$ elements in $A$, there exists a `small' set $E\subseteq B[X]$ such that, if a monic polynomial of degree two $P\in B[X]$ has the property that each $P(a_i)$ is a square in $B$, then $P\in E$ or $P$ is a square in $B[X]$.}

Also, we will show that such a result on representation of squares should hold in number fields (and hence B\"uchi's problem should also be true there) under the hypothesis that a  conjecture by Bombieri holds for surfaces\,:

\begin{conjecture}[Bombieri] \label{ConjBombieri} If $X$ is a smooth variety of general type defined over a number field $K/\Q$, then $X(K)$ is contained in a proper Zariski closed set of $X$.
\end{conjecture} 

Finally, we will use the positive answer obtained for B\"uchi's problem for $p$-adic meromorphic functions in order to improve some undecidability results for $p$-adic analytic and meromorphic functions. 

In the next section, we present in details the main results that are proven in this paper. 


\section{Main results}

\subsection{Representation of squares in the field of $p$-adic meromorphic functions}

Let $p$ be a prime number and let $\C_p$ be the field of $p$-adic complex numbers (the completion of the algebraic closure of the field $\Q_p$ of $p$-adic numbers). Throughout the paper, one can replace $\C_p$ by any algebraically closed field of characteristic zero, complete with respect to a non-trivial non-Archimedean valuation.

Let $\Acal_p$ be the ring of entire functions over $\C_p$ and let $\Mcal_p$ be the field of meromorphic functions over $\C_p$. We prove the following theorem on representation of squares by polynomials.

\begin{theorem}\label{RepMer} Let $P\in \Mcal_p[X]$ be a monic polynomial of degree two. If $P(a)$ is a square in $\Mcal_p$ for at least $35$ values of $a\in\C_p$, then either $P$ has constant coefficients or $P$ is a square in $\Mcal_p[X]$.
\end{theorem}

\begin{corollary}\label{BuchiMer} The problems $\BP(\Acal_p)$ and $\BP(\Mcal_p)$ have a positive answer. 
\end{corollary}

\noindent We will prove these results in Section \ref{SecMer}. 

\subsection{Representation of squares in number fields}

We will show in Sections \ref{geometric} and \ref{SecNum} the following results on representation of squares in number fields.

\begin{theorem}\label{RepNum} Assume Conjecture \ref{ConjBombieri} holds for surfaces. Let $K$ be a number field and $\{a_1,\ldots,a_8\}$ a set of eight elements in $K$. There exists a finite (possibly empty) set $E=E(K,(a_i)_i)$ of polynomials in $K[x]$ such that the following holds\,: for each polynomial $f$ of the form $x^2+ax+b\in K[x]$, if $f(a_i)$ are squares in $K$ for each $i$ then either $f\in E$, or $f=(x+c)^2$ for some $c\in K$.
\end{theorem}
It is an obvious but remarkable fact that if one could find a number field $K$ and a sequence $(a_i)$ such that the set $E(K,(a_i))$ is infinite, then one would automatically obtain a counterexample to Bombieri's Conjecture.

From the finiteness of the set $E(K,(a_i))$ one can easily derive the following (see Section \ref{SecNum}).

\begin{corollary}\label{buchistrong} Assume that Bombieri's conjecture holds for surfaces defined over $\Q$. Let $a_1,a_2,\ldots$ be a sequence of integers without repeated terms. There exists a constant $M$ (depending on the sequence $(a_i)_i$) such that: 
if a polynomial $f=x^2+ax+b\in \Q[x]$ satisfies the property `$f(a_i)$ is a square in $\Z$ for $i=1,\ldots,M$', then $f$ is of the form $f=(x+c)^2$, for some $c\in\Z$.
\end{corollary}

Observe that the dependence of $M$ on the sequence cannot be dropped. Consider for example the polynomial $f_N=x^2-4(2N)!$, and define 
$$
a_i=i! +\frac{(2N)!}{i!}.
$$ 
Then it is obvious that $(a_i)_{i=1}^N$ is a strictly decreasing sequence in $\Z$ and each $f_N(a_i)$ is a square in $\Z$. 

Note that, if in Corollary \ref{buchistrong} we set $a_n=n$ for each $n$, then we obtain a positive answer to B\"uchi's Problem for $\Z$ (under Bombieri's Conjecture).

For $\Z$, we can actually state a stronger result, which will be proved in Section \ref{SecNum}. In order to state it, we need first to give a definition. Observe that if a sequence $(x_n)$ satisfies the equation
$$
x_{n+2}^2-2x_{n+1}^2+x_n^2=2
$$
for $n=1,\dots,M-2$ then solving the induction gives
$$
x_n^2=(n-1)(n-2)-(n-2)x_1^2+(n-1)x_2^2
$$
for each $n$. This observation motivates the following:

\begin{definition} For $n\ge 3$, we say that a complex projective surface $X\subseteq\Pro^n$ is an \emph{$n$-B\"uchi surface}, if there exists a sequence of distinct nonzero \emph{integers} $\delta_2,\delta_3,\ldots,\delta_n$ such that $X$ is defined by the system 
$$
\delta_2 x^2_i=\delta_i\delta_2(\delta_i-\delta_2)x^2_0-(\delta_i-\delta_2)x^2_1+\delta_ix^2_2
$$
where $i$ ranges from $3$ to $n$ (the $x_i$ being the independent variables). 
\end{definition}

One can verify that all B\"uchi surfaces are smooth, and that $n$-B\"uchi surfaces for $n\ge 6$ are of general type (we will not prove this fact, because the proof is essentially the same as the one we give for Lemma \ref{Xn}). 

\begin{theorem}\label{LogicNum} If Bombieri's conjecture for $K=\Q$ holds for \emph{some} $n$-B\"uchi surface for $n\ge 8$, then multiplication is existentially definable in $\Z$ over $\Lcal_2$. 
\end{theorem}


Therefore, we need quite less than a positive answer to B\"uchi's problem in order to define multiplication (because B\"uchi's problem corresponds to the particular case where $\delta_k=k-1$ for each $k$).

\subsection{Representation of squares in function fields}

The geometric results in Section \ref{geometric} will be used in Section \ref{SecFun} to prove the following theorem, which is the analogue of Theorem \ref{RepMer}.

\begin{theorem}\label{RepFun} Let $F$ be a field of characteristic zero and $C$ a non-singular projective curve defined over $F$. Define the integer $M=\max\{8,4(g+1)\}$ where $g$ is the genus of $C$. Write $K(C)$ for the function field of $C$ and let $X$ be transcendental over $K(C)$. Let $P\in K(C)[X]$ be a monic polynomial of degree two. If $P(a)$ is a square in $K(C)$ for at least $M$ values of $a\in F$, then either $P$ has constant coefficients or $P$ is a square in $K(C)[X]$.
\end{theorem}

Theorem \ref{RepFun} gives as a direct consequence a positive answer to B\"uchi's problem, but such a positive answer is not new since it was (implicit) in \cite{Vojta2} and recently was solved by a new method in characteristic zero and (large enough) positive characteristic in \cite{ShlapentokhVidaux}.

\subsection{Undecidability for $p$-adic entire and meromorphic functions in B\"uchi's language}\label{ResLogic}

Corollary \ref{BuchiMer} allows us to obtain very strong undecidability results for $p$-adic analytic and meromorphic functions, improving results by Lipshitz, Pheidas and Vidaux. In order to state the theorems, we need to introduce some notation.

Recall that $\Acal_p$ stands for the ring of entire functions over $\C_p$ and $\Mcal_p$ stands for the field of meromorphic functions over $\C_p$, with variable $z$. 

Consider the language  $\Lcal_2'=\{0,1,+,f_z,P_2\}$ where $P_2(x)$ is interpreted as `$x$ is a square' and $f_z(x,y)$ is interpreted as `$y=zx$'. 
 
\begin{theorem}\label{LogicMer} Multiplication is positive existentially definable in $\Mcal_p$ and in $\Acal_p$ over the language $\Lcal_2'$.
\end{theorem}

See Section \ref{SecLogicMer} for a proof.

Define the languages $\Lcal_R^z=\{0,1,+,\cdot,z\}$, $\Lcal_R^*=\{0,1,+,\cdot,z,\ord\}$, $\Lcal_2^z=\{0,1,+,P_2,f_z\}$ and $\Lcal_2^*=\{0,1,+,P_2,f_z,\ord\}$. In the field $\Mcal_p$ and the ring $\Acal_p$, we interpret $\ord(x)$ as `$x(0)=0$'.

We recall that the following two theories are undecidable: the positive existential theory of $\Acal_p$ in the language $\Lcal_R^z$ (see \cite{LipshitzPheidas}) and the positive existential theory of $\Mcal_p$ in the language $\Lcal_R^*$ (see \cite{Vidaux3}). From this and Theorem \ref{LogicMer} we conclude

\begin{theorem} The positive existential theory of $\Acal_p$ in the language $\Lcal_2^z$ and the positive existential theory of $\Mcal_p$ in the language $\Lcal_2^*$ are undecidable.
\end{theorem}





\section{Some results in $p$-adic Nevanlinna Theory}

First we present the notation we use for the usual functions in Nevanlinna Theory. 

We will work over the field $\C_p$ with absolute value $|\cdot|_p$. Write $\Acal_p$ for the ring of entire functions over $\C_p$ and $\Mcal_p$ for the field of meromorphic functions over $\C_p$. We denote by $F^+$ the positive part of a function $F$ in to $\R$, that is $F^+=\max\{F,0\}$.
We adopt the following notation for the standard functions in $p$-adic Nevanlinna theory, where $f=\frac{h}{g}\in\Mcal_p$ is non-zero, and where $g,h\in\Acal_p$ are coprime:
\begin{eqnarray*}
B[r]&=&\{z\in\C_p\colon |z|_p\leq r\}\\
n(r,h,0)&=&\mbox{number of zeroes of }h\mbox{ in }B[r]\mbox{ counting multiplicity}\\
n(r,f,0)&=&n(r,h,0)\\
n(r,f,\infty)&=&n(r,g,0)\\
N(r,h,0)&=&\int_0^r\frac{n(t,h,0)-n(0,h,0)}{t}dt+n(0,h,0)\log r\\
N(r,f,0)&=&N(r,h,0)\\
N(r,f,a)&=&N(r,f-a,0)\\
N(r,f,\infty)&=& N(r,g,0)\\
|h|_r&=&\max_{n\ge 0}|a_n|_pr^n\\
|f|_r&=&\frac{|h|_r}{|g|_r}\\
m(r,f,a)&=&\log^+\frac{1}{|f-a|_r}\\
m(r,f,\infty)&=&\log^+|f|_r
\end{eqnarray*}

We recall to the reader that for each $r>0$, the function $|\cdot|_r:\Mcal\to \R$ is a non-archimedean absolute value satisfying $|a|_r=|a|_p$ when $a$ is constant.

We will need the following results from $p$-adic Nevanlinna Theory. For a general presentation of $p$-adic complex analysis, see for example \cite{Robert}.

First we have the \emph{Logarithmic Derivative Lemma} (see \cite{CherryYe}, Lemma 4.1):

\begin{lemma}\label{LDL} If $n>0$ is a positive integer and $f\in\Mcal_p$ then 
$$
\left|\frac{f^{(n)}}{f}\right|_r\le \frac{1}{r^n}
$$
where $f^{(n)}$ stands for the $n$-th derivative.
\end{lemma}

We will also need the \emph{Poisson-Jensen Formula} (see \cite{CherryYe}, Lemma 3.1):
\begin{theorem}\label{PJF} Given $f\in \Mcal_p$, there exists a constant $C$ depending only on $f$ such that
$$
\log|f|_r=N(r,f,0)-N(r,f,\infty)+C.
$$
\end{theorem}

As a consequence of the Poisson-Jensen Formula, we get the \emph{First Main Theorem}:

\begin{theorem}\label{FMT}
Let $f\in\Mcal_p$ be a non-constant meromorphic function and $a\in\C_p$. As $r\to\infty$ we have
$$m(r,f,a)+N(r,f,a)=m(r,f,\infty)+N(r,f,\infty)+\Ocal(1).$$
\end{theorem}

Finally, we state the \emph{Second Main Theorem} (see \cite{Ru}, Theorem 2.1):

\begin{theorem}\label{SMT} Let $f\in\Mcal_p$ be a non-constant meromorphic function and let $a_1,\ldots,a_q\in\C_p$ be distinct. Then, as $r\to\infty$ we have
$$
\sum_{i=1}^q m(r,f,a_i)\le N(r,f,\infty)+\Ocal(1).
$$
\end{theorem}


\section{Proof of Theorem \ref{RepMer} (Meromorphic Functions)} \label{SecMer}

The following equality will be used many times without mention within this section:
\begin{equation}\label{nN}
N(r,f,x)=K+\int_1^r\frac{n(t,f,x)}{t}dt,\quad \mbox{for large }r.
\end{equation}
It will be used systematically in order to deduce inequalities (for large $r$) about $N$ when we know inequalities about $n$ (the point is that the integral is a linear and monotone operator).\\

In order to simplify the proof of Theorem \ref{RepMer}, we actually will prove the following equivalent result.

\begin{theorem}\label{mer2} Let $h_1,\ldots,h_M$ be elements of $\Mcal_p$ such that at least one $h_i$ is non-constant. Let $a_1,\ldots,a_M$ be $M$ distinct elements of $\C_p$. If there exist $f,g\in\Mcal_p$, with $g$ non-zero, such that
\begin{equation}\label{Pizarra}
h_j^2=(a_j+f)^2-g \qquad j=1,\ldots,M
\end{equation}
then $M\le 34$.
\end{theorem}

For the rest of this section, we will assume that we are under the hypothesis of Theorem \ref{mer2}. Assuming $M\ge 35$ we will obtain a contradiction.\\

First, we observe that
\begin{eqnarray}\label{diferencia}
h_i^2-h_j^2=(a_i-a_j)(2f+ a_i+a_j).
\end{eqnarray}

\begin{lemma}\label{fnocte} The function $f$ is not constant.
\end{lemma}
\begin{proof} If $f$ is constant then so is $c_i=(a_i+f)^2$. Note that since some $h_i$ is non-constant, $g$ is non-constant. Taking $i$, $j$ and $k$ such that $c_i$, $c_j$, and $c_k$ are pairwise distinct constants, the following equality
$$
(h_ih_jh_k)^2=(c_i-g)(c_j-g)(c_k-g)
$$
gives a non-constant meromorphic parametrization of an elliptic curve over $\C_p$, which is impossible by a theorem of Berkovich (see \cite{Berkovich}).
\end{proof}

\begin{lemma}\label{poles} Let $x\in \C_p$ be a pole of some $h_i$. There exists an index $k$ depending on $x$ such that for each $i\ne k$ we have (simultaneously)
\begin{enumerate}
\item $\ord_x h_k\ge \ord_x h_i$; 
\item $\ord_x f\ge 2\ord_x h_i$;
\item $\ord_x g\ge 4\ord_x h_i$; 
\item $\ord_x h_i=\ord_x h_j$ for all $j\ne k$; and
\item $\ord_x h_i\le -1$.
\end{enumerate}
Moreover, for each $i$ we have 
\begin{equation}\label{polesord}
\min\{\ord_x h_i,0\}\geq\frac{1}{M-1}\sum_l \min\{\ord_x h_l,0\}
\end{equation}
and, there exists a positive constant $K$ such that for large enough $r$ and for each $i$ we have
\begin{equation}\label{polesN}
N(r,h_i,\infty)\le \frac{1}{M-1}\sum_l N(r,h_l,\infty).
\end{equation}
\end{lemma}
\begin{proof}
Let $i_0$ be an index such that $h_{i_0}$ has a pole at $x$. 

First suppose that \emph{all} $h_i$ have the same order at $x$ (hence negative). In this case, Items (1), (4) and (5) hold trivially, Item (2) comes from Equation \eqref{diferencia}, and Item (3) comes from Equation \eqref{Pizarra}. Indeed for Item (3) we have
$$
\begin{aligned}
\ord_x(g)&\geq2\min\{\ord_x(h_i),\ord_x(f+a_i)\}\\
&=2\min\{\ord_x(h_i),\ord_x(f)\}\\
&=2\min\{\ord_x(h_i),2\ord_x (h_i)\}\\
&\geq4\ord_x h_i,
\end{aligned}
$$
where the last inequality comes from Item (2). 

The other case is when not all $h_i$ have the same order at $x$. Choose $k$ such that item (1) holds true. By Equation \eqref{diferencia} for indices $k$ and any $i\ne k$, Item (4) holds true. If $i_0=k$ then all $h_i$ have a pole at $x$ (by maximality of $k$), and if $i_0\ne k$ then by Item (4), for all $i\ne k$, $h_i$ has a pole at $x$. Hence Item (5) holds true. Items (2) and (3) for $i\ne k$ follow as in the previous case.

Finally, by Items (1), (4) and (5), and observing that $\ord_x h_k$ could be positive, we have for each $i$ 
$$
(M-1)\min\{\ord_x h_i,0\}=\sum_{l\ne k}\min\{\ord_x h_l,0\}\geq\sum_l \min\{\ord_x h_l,0\}.
$$
Summing for $x\in B[r]$ we obtain
$$
(M-1)n(r,h_i,\infty)\leq \sum_l n(r,h_l,\infty).
$$
which gives the inequality \eqref{polesN} by Equation \eqref{nN}.
\end{proof}

\begin{lemma}\label{asterisco} The following inequality holds
$$
\sum_{n=1}^{M} \log|h_n|_r +\frac{1}{M-1}\sum_{n=1}^{M} N(r,h_n,\infty) \ge -\frac{1}{2}N(r,f,\infty)+\Ocal(1).
$$
\end{lemma}
\begin{proof} By the Second Main Theorem \ref{SMT}, we have for each $i\in\{1,\ldots,M\}$
\begin{eqnarray*}
-N(r,f,\infty)+\Ocal(1) \le -\sum_{j\ne i}\log^+\left|\frac{1}{f+\frac{a_i+a_j}{2}}\right|_r
\le \sum_{j\ne i}\log\left|f+\frac{a_i+a_j}{2}\right|_r.
\end{eqnarray*}
Since by Equation \eqref{diferencia} we have 
$$
h_i^2-h_j^2=2(a_i-a_j)\left(f+\frac{a_i+a_j}{2}\right),
$$ 
we deduce
$$-N(r,f,\infty)+\Ocal(1) \le \sum_{j\ne i}\log\left|h_i^2-h_j^2\right|_r.$$
If for a given $r$, $i_r$ is an index such that $|h_i|_r$ is minimal, then 
\begin{eqnarray*}
\frac{1}{2}\sum_{j\ne i_r}\log\left|h_{i_r}^2-h_j^2\right|_r&\le
& \sum_{j\ne i_r}\log\left|h_j\right|_r \\
&=&C+ \sum_{j\ne i_r} \left(N(r,h_j,0)-N(r,h_j,\infty)\right)\\
&\le& C+N(r,h_{i_r},\infty)+ \sum_n \left(N(r,h_n,0)-N(r,h_n,\infty)\right)\\
&=&C'+N(r,h_{i_r},\infty)+ \sum_n \log\left|h_n\right|_r\\
&\le& C''+\frac{1}{M-1}\sum_n N(r,h_n,\infty)+ \sum_n \log\left|h_n\right|_r
\end{eqnarray*}
where the first and second equalities are given by the Poisson-Jensen Formula \ref{PJF}, the third inequality is given by Lemma \ref{poles} (see Equation \eqref{polesN}), and $C$, $C'$, $C''$ are fixed constants (not depending on $r$ nor on $i_r$).

Finally we have
$$
-\frac{1}{2}N(r,f,\infty)+\Ocal(1) \le \frac{1}{2}\sum_{j\ne i_r}\log\left|h_{i_r}^2-h_j^2\right|_r\le \sum \log\left|h_n\right|_r + \frac{1}{M-1}\sum N(r,h_n,\infty) + C''
$$
for each $r$ large enough, and the lemma is proven.
\end{proof}
\begin{lemma}\label{L} The following inequalities hold:
$$
n(r,f,\infty)\le \frac{2}{M-1}\sum_n n(r,h_n,\infty)
$$
and
$$
\sum_n N(r,h_n,0)\ge \frac{M-3}{M-1}\sum_n N(r,h_n,\infty)+\Ocal(1).
$$
\end{lemma}
\begin{proof} Observe that by Lemma \ref{poles} (Item (2) and Equation \eqref{polesord}) we have 
$$
(M-1)n(r,f,\infty)\le 2\sum n(r,h_j,\infty),
$$
hence
$$
(M-1)N(r,f,\infty)\le 2\sum N(r,h_n,\infty)+\Ocal(1).
$$
The second formula comes immediately by Lemma \ref{asterisco} and the Poisson-Jensen Formula \ref{PJF}.
\end{proof}
\\

The equations
\begin{eqnarray*}
h_n^2+g&=&(a_n+f)^2\\
2h'_nh_n+g'&=&2f'(a_n+f)
\end{eqnarray*}
are directly deduced by reordering and differentiating the one given in the hypothesis. From this we deduce
\begin{eqnarray*}
(2h'_nh_n+g')^2=4f'^2(h_n^2+g)
\end{eqnarray*}
hence
\begin{eqnarray*}
g'^2-4f'^2g=4h_n(h_nf'^2-h'^2_nh_n-h_n'g').
\end{eqnarray*}
Writing
\begin{eqnarray*}
\Delta&=&g'^2-4f'^2g\\
\Delta_n&=&h_nf'^2-h'^2_nh_n-h'_ng'
\end{eqnarray*}
we have
\begin{equation}\label{Deltas}
\Delta=4h_n\Delta_n.
\end{equation}
\begin{lemma}\label{C1} If $\Delta$ is not identically zero, then
$$N(r,\Delta,0)\ge \frac{1}{2}\sum N(r,h_n,0)-\frac{8}{M-1}\sum N(r,h_n,\infty)+\Ocal(1).$$ 
\end{lemma}
\begin{proof} On the one hand, for a given $x\in\C_p$ suppose $f$ has a pole at $x$ and $h_j(x)=0$ for some index $j$. Set $l=\ord_x (h_j)$ and $m=\ord_x(f)$. Note that $\ord_x(g)=2m$ because $h_j(x)=0$ (see Equation \eqref{Pizarra}). Write 
$$
h_j=u_l(z-x)^l+u_{l+1}(z-x)^{l+1}+\cdots,
$$ 
$$
f=v_m(z-x)^m+v_{m+1}(z-x)^{m+1}+\cdots
$$ 
and 
$$
g=w_{2m}(z-x)^{2m}+w_{2m+1}(z-x)^{2m+1}+\cdots
$$ 
for the Laurent series of $h_j$, $f$ and $g$ at $x$. Observe that $w_{2m}=v^2_m$. The first term of the Laurent series at $x$ for respectively $h_jf'^2$, $h'^2_jh_j$ and $h'_jg'$ is, respectively,  
\begin{eqnarray*}
m^2u_lv^2_m(z-x)^{l+2m-2}\\
l^2u_l^3(z-x)^{3l-2}\\
2lm u_lv^2_{m}(z-x)^{l+2m-2}
\end{eqnarray*}
hence $\ord_x\Delta_j=l+2m-2$ since $2l\ne m$. Therefore, we have 
$$
\ord_x\Delta=2(l+m-1).
$$

On the other hand, if $x\in\C_p$ is not a pole of $f$ and is a zero of some $h_j$, then we have 
$$
\ord_x\Delta\geq\ord_x(h_j)
$$ 
because by Equation \eqref{Pizarra}, $g$ does not have a pole, hence $\Delta_j$ does not have a pole and we conclude by Equation \eqref{Deltas}.

Let $A_r$ be the set of $x\in B[r]$ such that $f$ has not a pole at $x$ and $h_j(x)=0$ for some index $j$, and let $B_r$ be the set of $x\in B[r]$ such that $f$ has a pole at $x$ and $h_j(x)=0$ for some index $j$. Observe that, by Equation \eqref{diferencia}, no three of the $h_n$ can share a zero (we use it for the fifth inequality below). We have then
\begin{eqnarray*}
n(r,\Delta,0)&\ge& \sum_{x\in A_r}\ord_x\Delta + \sum_{x\in B_r}\ord_x\Delta\\
&\ge& \sum_{x\in A_r}\max_{h_i(x)=0} \ord_x(h_i) + \sum_{x\in B_r}\max_{h_i(x)=0} 2(\ord_x(h_i)+\ord_x(f)-1)\\
&\ge& \sum_{x\in A_r\cup B_r}\max_{h_i(x)=0} \ord_x(h_i) + 2\sum_{x\in B_r} \max_{h_i(x)=0}(\ord_x(f)-1)\\
&=& \sum_{x\in A_r\cup B_r}\max_{h_i(x)=0} \ord_x(h_i) + 2\sum_{x\in B_r} (\ord_x(f)-1)\\
&\ge& \sum_{x\in A_r\cup B_r}\max_{h_i(x)=0} \ord_x(h_i) + 4\sum_{x\in B_r} \ord_x(f)\\
&\ge& \frac{1}{2}\sum_i n(r,h_i,0) - 4n(r,f,\infty)\\
&\ge& \frac{1}{2}\sum_i n(r,h_i,0) - \frac{8}{M-1}\sum n(r,h_i,\infty)
\end{eqnarray*}
where the last inequality comes from Lemma \ref{L}. The result follows.
\end{proof}



\begin{lemma}\label{C2} If $\Delta$ is not identically zero, then
$$N(r,\Delta,\infty)\le \frac{8}{M-1}\sum N(r,h_n,\infty)+\Ocal(1).$$
\end{lemma}
\begin{proof} Suppose that some $x\in\C_p$ is a pole of $\Delta$. Then, by definition of $\Delta$, it is a pole of  $f$ or of $g$. If none of the $h_i$ has a pole at $x$ then by Equation \eqref{diferencia} $f$ does not have a pole, and by Equation \eqref{Pizarra}, $g$ does not have a pole, which contradicts our hypothesis. Therefore, some $h_i$ has a pole at $x$. Take $k$ as in Lemma \ref{poles}. For each index $i\ne k$ we have (observing that $\ord_x(h_i)\leq -1$ and that if $g'=0$ then $\ord_x h_i'g'$ is infinite)
\begin{eqnarray*}
\ord_x\Delta &\ge& \ord_x h_i + \min\{\ord_x h_if'^2,\ord_x h'^2_ih_i,\ord_x h'_ig'\}\\
&\ge& \ord_x h_i +\min\{7\ord_x h_i,5\ord_x h_i,7\ord_x h_i\}\\
&=&8\ord_x h_i.
\end{eqnarray*}
Hence, using the Lemma \ref{poles} (Equation \eqref{polesord}) we have
$$
\ord_x\Delta\ge \frac{8}{M-1}\sum_l \min\{\ord_x h_l,0\}.
$$

Write $D_r$ for the set of poles of $\Delta$ in $B[r]$. We have
\begin{eqnarray*}
n(r,\Delta,\infty)&=&\sum_{x\in D_r}-\ord_x\Delta \\
&\le& \frac{8}{M-1}\sum_{x\in D_r}\sum_l \max\{-\ord_x h_l,0\}\\
&\le& \frac{8}{M-1}\sum_l \sum_{x\in B[r]} \max\{-\ord_x h_l,0\}\\
&=&\frac{8}{M-1}\sum_l n(r,h_l,\infty).
\end{eqnarray*}
and the result follows.
\end{proof}
\begin{lemma}\label{PeloCorto}
\begin{enumerate}
\item For each $r>0$, there exists an index $k_r$ such that $|h_{k_r}|_r$ is minimal. 
\item There exists a positive constant $K_f$ such that, for any $r>0$ and for all $i\ne k_r$, we have
$$
\log |f|_r\le \max\{0,2\log |h_i|_r\} + K_f.
$$
\item There exists a positive constant $K_g$ such that, for any $r>0$ and for all $i\ne k_r$, we have
$$
\log |f|_r\le \max\{0,4\log |h_i|_r\} + K_g.
$$
\end{enumerate}
\end{lemma}
\begin{proof}
Item (1) is immediate since for each $r$, the set $\{|h_i|_r\colon i=1,\dots,M\}$ is finite. Let us prove Item (2). There exists a positive constant $K'>1$ such that for each $r>0$, $i$ and $j$, we have
\begin{equation}\label{Pelado}
|2f|_r\le |2f+a_i+a_j|_r+|a_i+a_j|_r\leq K'+ |2f+a_i+a_j|_r.
\end{equation}
On the other hand, by Equation \eqref{diferencia} there exists a constant $K''>1$ such that, for any $r>0$, $i\ne k_r$ and $j$, we have
$$
\begin{aligned}
|2f+a_i+a_j|_r&=\left|\frac{h_i^2-h_{k_r}^2}{a_i-a_{k_r}}\right|_r&\qquad\\
&\leq \left|\frac{h_i^2}{a_i-a_{k_r}}\right|_r &\qquad \textrm{(by Item (1))}\\
&\leq  K''|h_i^2|_r&
\end{aligned}
$$
hence by Equation \eqref{Pelado}
$$
|2f|_r\leq K''|h_i^2|_r+K'\leq K'' \max\{|h_i^2|_r,1\}+ K'.
$$
Therefore, we have
$$
\begin{aligned}
\log|f|_r&\leq \log(K'' \max\{|h_i^2|_r,1\}+ K')-\log|2|_r\\
&\leq \log(K'' \max\{|h_i^2|_r,1\})+\log K'+\log2-\log|2|_r\\
&\leq \max\{2\log |h_i|_r,0\} + K_f
\end{aligned}
$$
with $K_f$ is a positive constant greater than $\log K''+\log K'+\log2-\log|2|_r$, and where the second inequality comes from the fact that for all real numbers $x,y\geq 1$, we have $\log(x+y)\leq\log x+\log y+\log2$ (just write $(1-x)(y-1)\le 0$).

Finally, we prove Item (3). By Equation \eqref{Pizarra} and Item (2), for each $i\ne k_r$ we have
$$
\begin{aligned}
\log |g|_r&=\log|(f+a_i)^2-h_i^2|_r \\
&\le \log\left(\max\{|h_i^2|_r, |f^2|_r,|2a_if|_r,|a_i^2|_r\}\right)\\
&\le \log\left(\max\{|h_i^2|_r, |f^2|_r,|a_i^2|_r\}\right)\\
&\le \max\{2\log |h_i|_r, 0, 2\log |f|_r\} + \max\{|a_i^2|_r\} \\
&\le \max\{2\log |h_i|_r, 0, 2\max\{0,2\log |h_i|_r\} + 2K_f\} + \max\{|a_i^2|_r\} \\
&\le \max\{2\log |h_i|_r+2K_f,2K_f, 4\log |h_i|_r+2K_f\} + \max\{|a_i^2|_r\} \\
&\le \max\{0,4\log |h_i|_r\} + K_g
\end{aligned}
$$
where $K_g$ is a fixed positive constant bigger than $2K_f+\max\{|a_i^2|_r\}$, and where the second inequality comes from the following\,:
$$
|2a_if|_r\le|a_i|_r|f|_r\leq\frac{|a_i^2|_r+|f^2|_r}{2}\leq\max\{|a_i^2|_r,|f^2|_r\}.
$$

\end{proof}

\begin{lemma}\label{C3} If $\Delta$ is not identically zero, then 
$$\log |\Delta|_r\le \frac{6M-2}{M(M-1)}\sum \log |h_n|_r+\frac{8}{(M-1)^2}\sum N(r,h_n,\infty)-2\log r+\Ocal(1).$$
\end{lemma}
\begin{proof} By the Poisson-Jensen Formula \ref{PJF} and Lemma \ref{poles} (Equation \eqref{polesN}) we have for $r$ large enough and for each $i$
$$
\begin{aligned}
\log|h_{i}|_r&= N(r,h_{i},0)-N(r,h_{i},\infty)+C\\
&\ge -N(r,h_{i},\infty)+C\\
&\ge -\frac{1}{M-1}\sum_n N(r,h_n,\infty)+C'
\end{aligned}
$$
for some constants $C,C'$. 

Given $r>0$ take $k_r$ as in Lemma \ref{PeloCorto}. Choose $i_r$ such that $|h_{i_r}|_r$ is minimal in $\{|h_i|_r\colon i\ne k_r\}$, and note that
$$
\log |h_{i_r}|_r\le \frac{1}{M-1}\sum_{i\ne k_r}\log |h_i|_r.
$$

By Item (2) in Lemma \ref{PeloCorto}, we have for each $r$ large enough
$$
\begin{aligned} 
\log |f|_r&\le \max\{0,2\log |h_{i_r}|_r\} + K_f\\
&\le \max\left\{0,\frac{2}{M-1}\sum_{i\ne k_r}\log |h_i|_r\right\} + K_f\\
&\le \max\left\{0,\frac{2}{M-1}\sum \log\left|h_n\right|_r + \frac{2}{(M-1)^2}\sum N(r,h_n,\infty)\right\} - \frac{2C'}{M-1} + K_f
\end{aligned}
$$
and by Item (3) in Lemma \ref{PeloCorto} we have for each $r$ large enough
$$
\begin{aligned}
\log |g|_r&\le \max\{0,4\log |h_{i_r}|_r\} + K_g\\
&\le \max\left\{0,\frac{4}{M-1}\sum_n \log\left|h_n\right|_r + \frac{4}{(M-1)^2}\sum_n N(r,h_n,\infty)\right\} - \frac{4C'}{M-1}+K_g.
\end{aligned}
$$

Hence, for large $r$ we get
\begin{eqnarray} 
\label{cotaf}\log |f|_r &\le& \max\left\{0,\frac{2}{M-1}\sum_n \log\left|h_n\right|_r + \frac{2}{(M-1)^2}\sum_n N(r,h_n,\infty)\right\} + \Ocal(1)\\
\label{cotag}\log |g|_r &\le& \max\left\{0,\frac{4}{M-1}\sum_n \log\left|h_n\right|_r + \frac{4}{(M-1)^2}\sum_n N(r,h_n,\infty)\right\} + \Ocal(1).\\
\end{eqnarray}

Since $\Delta$ is not the zero function, from Lemma \ref{LDL} (Logarithmic Derivative Lemma) we have for each index $n$
\begin{eqnarray*}
|\Delta|_r\le |h_n|_r\max\{|h_nf'^2|_r,|h'^2_nh_n|_r,|h'_ng'|_r\}\le \frac{1}{r^2}|h_n|^2_r\max\{|f|^2_r,|h_n|^2_r,|g|_r\}
\end{eqnarray*}
and by Equation \eqref{Pizarra} for each $n$ holds 
$$
2\log |h_n|_r\le \max\{2\log |f|_r, 0 , \log |g|_r\}+\Ocal(1)
$$
therefore we have for each $n$
\begin{eqnarray*}
\log |\Delta|_r\le \log\left(\frac{1}{r^2}|h_n|^2_r\right)+ \max\{2\log |f|_r, 0 , \log |g|_r\}+\Ocal(1).
\end{eqnarray*}
Since this last expression is true \emph{for each} $n$, we have
\begin{eqnarray*}
\log |\Delta|_r\le \frac{2}{M}\sum\log |h_n|_r - 2\log r+ \max\{2\log |f|_r, 0 , \log |g|_r\}+\Ocal(1).
\end{eqnarray*}

Note that by equations \eqref{cotaf} and \eqref{cotag} 
$$
\max\{2\log |f|_r, 0 , \log |g|_r\}\le \max\left\{0,\frac{4}{M-1}\sum \log\left|h_n\right|_r + \frac{4}{(M-1)^2}\sum N(r,h_n,\infty)\right\} + \Ocal(1)
$$
and by Lemma \ref{asterisco} we have that this last expression is less than or equal to
$$
\frac{4}{M-1}\sum \log\left|h_n\right|_r + \frac{4}{(M-1)^2}\sum N(r,h_n,\infty)+\frac{2}{M-1}N(r,f,\infty)+\Ocal(1).
$$

Therefore
\begin{eqnarray*}
\log |\Delta|_r\le \left(\frac{2}{M}+\frac{4}{M-1}\right)\sum\log |h_n|_r - 2\log r + \frac{4}{(M-1)^2}\sum N(r,h_n,\infty)\\
 + \frac{2}{M-1}N(r,f,\infty)+\Ocal(1)
\end{eqnarray*}
Finally, we bound $N(r,f,\infty)$ using Lemma \ref{L} and the result follows.
\end{proof}
\begin{lemma}\label{eqdif} $\Delta=0$ for $M\ge 35$.
\end{lemma}
\begin{proof} The proof goes by contradiction, so we assume $\Delta$ is not identically zero. Consider the equation 
$$
\log|\Delta|_r=N(r,\Delta,0)-N(r,\Delta,\infty) +\Ocal(1).
$$
Lemmas \ref{C1}, \ref{C2} and \ref{C3} allow us to bound $\log |\Delta|_r$ above and below, obtaining
\begin{eqnarray*}
\frac{6M-2}{M(M-1)}\sum \log |h_n|_r+\frac{8}{(M-1)^2}I-2\log r \ge \frac{1}{2}Z-\frac{8}{M-1}I - \frac{8}{M-1}I +\Ocal(1) 
\end{eqnarray*}
where we write $Z=\sum N(r,h_n,0)$ and $I=\sum N(r,h_n,\infty)$. Observe that 
$$
\sum \log |h_n|_r= Z-I+\Ocal(1)$$
by Poisson-Jensen Formula \ref{PJF}. This and Lemma \ref{L} give
\begin{eqnarray*}
-2\log r &\ge& \left(\frac{1}{2}-\frac{6M-2}{M(M-1)}\right)Z +\left(\frac{6M-2}{M(M-1)} -\frac{16}{M-1} -\frac{8}{(M-1)^2}\right)I+\Ocal(1)\\
&\ge& \left(\left(\frac{1}{2}-\frac{6M-2}{M(M-1)}\right)\frac{M-3}{M-1} -\frac{10M^2-2}{M(M-1)^2} \right)I+\Ocal(1)\\
&=&\frac{M^2-35M+8}{2M(M-1)}I+\Ocal(1).
\end{eqnarray*}
A contradiction for $M\ge 35$.
\end{proof}

Finally, we have that $\Delta$ is the zero function, $f$ is not a constant and $g$ is non-zero. We get the equation $g'^2=4f'^2g$, which implies that exists a meromorphic function $u$ such that $g=u^2$ and $u'^2=f'^2$. Hence $u= \alpha f + b$ for certain $\alpha\in\{-1,1\}$ and $b\in\C_p$, and we obtain 
\begin{eqnarray*}
h_n^2&=&(a_n+f)^2-(\alpha f+b)^2\\
&=&(a_n+f)^2-(f+\alpha b)^2\\
&=&(a_n-\alpha b)(a_n+\alpha b+2f).
\end{eqnarray*}
From this we get
\begin{eqnarray*}
\left(\frac{h_ih_jh_k}{\sqrt{(a_i-\alpha b)(a_j-\alpha b)(a_k-\alpha b)}}\right)^2=(a_i+\alpha b +2f)(a_j+\alpha b +2f)(a_k+\alpha b +2f).
\end{eqnarray*}
This is a contradiction because $f$ is not a constant. Therefore the Theorem \ref{mer2} is proven.


\section{Some geometric results}\label{geometric}

This section contains most of the geometric results that we will use in the next two sections. The arguments given here essentially are adaptations of the arguments given by Vojta in \cite{Vojta2}. Anyway, we prefer to perform most of the computations in order to give a clear presentation.

During the whole section, we assume that the base field is $\C$, and we write $g(X)$ for the genus of the curve $X$.

Let $S=(\delta_2,\delta_3,\ldots)$ be a sequence in $\C^*$ with pairwise distinct terms. Set $X_2=\Pro^2(\C)$ and for $n>2$ let $X_n\subset\Pro^n(\C)$ be the algebraic set defined by the equations
\begin{eqnarray}\label{defining}
\delta_2 x^2_i=\delta_i\delta_2(\delta_i-\delta_2)x^2_0-(\delta_i-\delta_2)x^2_1+\delta_ix^2_2
\end{eqnarray}
as $i$ ranges from $3$ to $n$. If $[x_0:\cdots:x_n]\in X_n$ is easy to see that at most $2$ of the $x_i$ can be zero, hence $X_n\subset U_0\cup U_1\cup U_2$ where $U_i$ is the open set $\{x_i\ne 0\}$. 
\begin{lemma}\label{Xn} $X_n$ is a smooth surface in $\Pro^n$, contains the lines
\begin{eqnarray}\label{triviallines}
\pm x_1=\pm x_2 -\delta_2 x_0=\cdots=\pm x_n - \delta_n x_0
\end{eqnarray}
and has canonical sheaf $\Ocal_{X_n}(n-5)$. In particular, $X_n$ is of general type for $n\ge 6$.
\end{lemma}
\begin{proof} Observe that, for $[x_0:\cdots,x_n]\in X_n\cap U_0$ the matrix 
\begin{eqnarray}\label{matrix}
\left[\begin{array}{ccccccc}
(\delta_3-\delta_2)x_1 &-\delta_3 x_2 &\delta_2 x_3 &0            &\cdots     &0           \\
(\delta_4-\delta_2)x_1 &-\delta_4 x_2 &0            &\delta_2 x_4 &\ddots     &0           \\
\vdots                 &\vdots        &\vdots       &\ddots       &\ddots     &\vdots      \\
(\delta_n-\delta_2)x_1 &-\delta_n x_2 &0            &0            &\cdots     &\delta_2 x_n 
\end{array}\right]
\end{eqnarray}
has rank $n-2$; indeed, we have $3$ cases depending on the number of zeroes between $x_3,\ldots, x_n$:
\begin{enumerate}
\item No zero: trivial.
\item One zero: at least one of the first two columns has no zero.
\item Two zeroes: suppose that $x_i=x_j=0$ where $3\le i<j\le n$, then no entry in the first two columns is zero. Therefore
\begin{eqnarray*}
\left|\begin{array}{cc}
(\delta_i-\delta_2)x_1 &-\delta_i x_2 \\
(\delta_j-\delta_2)x_1 &-\delta_j x_2
\end{array}\right|=\delta_2x_1x_2(\delta_j-\delta_i)\ne 0.
\end{eqnarray*}
\end{enumerate}
hence, $X_n$ is nonsingular at each point in $X_n\cap U_0$. The verification that $X_n$ is nonsingular at each point in $X_n\cap U_1$ and $X_n\cap U_2$ is quite similar, but the determinants at case (3) are 
\begin{eqnarray*}
\left|\begin{array}{cc}
-\delta_i\delta_2(\delta_i-\delta_2)x_0 &-\delta_i x_2 \\
-\delta_j\delta_2(\delta_j-\delta_2)x_0 &-\delta_j x_2
\end{array}\right|=\delta_2\delta_i\delta_jx_0x_2(\delta_j-\delta_i)\ne 0
\end{eqnarray*} 
and 
\begin{eqnarray*}
\left|\begin{array}{cc}
-\delta_i\delta_2(\delta_i-\delta_2)x_0 &(\delta_i-\delta_2)x_1 \\
-\delta_j\delta_2(\delta_j-\delta_2)x_0 &(\delta_j-\delta_2)x_1
\end{array}\right|=\delta_2x_0x_1(\delta_j-\delta_i)(\delta_j-\delta_2)(\delta_i-\delta_2)\ne 0
\end{eqnarray*}
respectively. Therefore $X_n$ is a smooth surface in $\Pro^n$.\\
The claim about the lines \eqref{triviallines} is an easy computation (look at $U_0\cap X_n$).\\
Finally, since $X_n$ is a complete intersection surface in $\Pro^n$ defined as the intersection of $n-2$ smooth hypersurfaces of degree 2, its canonical sheaf is $\Ocal(2(n-2)-n-1)=\Ocal(n-5)$. 
\end{proof} 
\begin{definition} Define the trivial lines of $X_n$ as the lines \eqref{triviallines}.
\end{definition}
Observe that for $n\ge 3$ the rational map $[x_0:\cdots:x_{n}]\mapsto [x_0:\cdots:x_{n-1}]$ induces a finite morphism $\pi_n:X_n\rightarrow X_{n-1}$ of degree $2$ ramified along the curve $C_n\subset X_n$ defined by $x_n=0$. This curve is nonsingular; indeed, if $[x_0:\cdots:x_n]\in C_n=X_n\cap\{x_n=0\}$ then at most one of the $x_0,\ldots,x_{n-1}$ can be zero and the remaining verification can be performed as in the proof of Lemma \ref{Xn} for cases (2) and (3) since $x_n=0$, but adding the extra row $(0,\ldots,0,1)$ to each matrix.\\
Define $\phi_n=\pi_3\circ\cdots\circ\pi_n$, we note that the image of $C_n$ in $X_2$ via $\phi_n$ is
\begin{eqnarray}\label{CninP2}
\delta_n\delta_2(\delta_n-\delta_2)x^2_0-(\delta_n-\delta_2)x^2_1+\delta_n x^2_2=0
\end{eqnarray}
\begin{definition} Let $X$ be a smooth surface over $\C$ and let $\Lcal$ be an invertible sheaf on $X$. Take a section $\omega\in H^0(X,\Lcal\otimes S^2(\Omega^1_X))$. Let $Y\subset X$ be a curve with normalization $i:\tilde{Y}\rightarrow Y$. We say that $Y$ is \emph{$\omega-$integral} if $i^*\omega\in H^0(\tilde{Y},i^*(\Lcal)\otimes S^2(\Omega^1_{\tilde{Y}}))$ vanishes identically on $\tilde{Y}$.
\end{definition}
On $U_0\subset\Pro^2=X_2$ define 
\begin{eqnarray*}
\omega=x_1x_2dx_1\otimes dx_1+(\delta^2_2-x_1^2-x_2^2)dx_1\otimes dx_2+ x_1x_2dx_2\otimes dx_2
\end{eqnarray*}
Note that, after the change of variables $y_0=x_0/x_1,y_2=x_2/x_1$, on $U_0\cap U_1$ we have
\begin{eqnarray*}
\omega=\frac{1}{y_0^5}\left(\delta_2^2y_0y_2dy_0\otimes dy_0+(1-\delta^2_2y_0^2-y_2^2)dy_0\otimes dy_2+ y_0y_2dy_2\otimes dy_2\right)
\end{eqnarray*}
hence $\omega$ extends to a section 
$$
\omega_2\in H^0(X_2,\Ocal_{X_{2}}(5)\otimes S^2(\Omega^1_{X_2})).
$$ 
\begin{lemma}\label{2integral} Write $[x_0:x_1:x_2]$ for homogeneous coordinates on $\Pro^2=X_2$. The only $\omega_2-$integral curves on $X_2$ are
\begin{enumerate}
\item $x_0=0$, $x_1=0$, and $x_2=0$
\item the four trivial lines
\item the conics $\delta_2c(c-\delta_2)x^2_0 - (c-\delta_2)x^2_1 + cx^2_2=0$ for $c\ne 0,\delta_2$.
\end{enumerate}
\end{lemma} 
\begin{proof} It is easy to see that curves of type (1) and (2) are $\omega_2-$integral. Let's show that curves of type (3) are $\omega_2$-integral. If we look at the affine chart $U_0$, on a curve of type (3) we have
$$
(c-\delta_2)x_1dx_1=cx_2dx_2
$$
hence
\begin{eqnarray*}
\omega_2&=&\left( \frac{c^2x^3_2}{(c-\delta_2)^2x_1} + \frac{cx_2}{(c-\delta_2)x_1}(\delta^2_2-x^2_1-x^2_2) + x_1x_2\right)dx_2\otimes dx_2\\
&=&\left( c^2x^2_2 + c(c-\delta_2)(\delta^2_2-x^2_1-x^2_2) + (c-\delta_2)^2x_1^2\right)\frac{x_2dx_2\otimes dx_2}{(c-\delta_2)^2x_1}\\
&=&\left( \delta^2_2c(c-\delta_2) - \delta_2(c-\delta_2)x^2_1 + \delta_2cx^2_2  \right)\frac{x_2dx_2\otimes dx_2}{(c-\delta_2)^2x_1}\\
&=&\delta_2\left( \delta_2c(c-\delta_2) - (c-\delta_2)x^2_1 + cx^2_2  \right)\frac{x_2dx_2\otimes dx_2}{(c-\delta_2)^2x_1}=0.
\end{eqnarray*}

Conversely, let $Y$ a $\omega_2-$integral curve on $X_2$ not of type (1) or (2), we will show that $Y$ is of type (3). Let $P\in Y$ be a regular point of $Y$ not in a line of type (1) nor (2). As $Y$ is regular at $P$, in some neighborhood of $P$ one can assume that one affine coordinate is function of the other, say $x_1=x_1(x_2)$. Since $Y$ is $\omega_2-$integral, we get a quadratic ordinary differential equation for $x_1$, hence there are $2$ local solutions at $P$. But exactly $2$ curves of type (3) pass through $P$. Therefore $Y$ is locally of type (3) on a dense set of points, so $Y$ is of type (3).
\end{proof}

Observe that the image of $C_n$ in $X_2$ is $\omega_2-$integral (see Equation \eqref{CninP2}).\\
Write $\omega'_n=\phi_n^*\omega_2$ and note that 
$$
\omega'_n\in H^0(X_n,\Ocal_{X_{n}}(5)\otimes S^2(\Omega^1_{X_n}))
$$
because $\pi_n^*\Ocal_{X_{n-1}}(1)=\Ocal_{X_n}(1)$ for each $n\ge 3$. 
\begin{lemma}\label{nintegral} Let $n\ge 6$ be an integer. The only $\omega'_n-$integral curves on $X_n$ are
\begin{enumerate}
\item the pull-backs via $\phi_n$ of the coordinate axes on $X_2$ to $X_n$
\item the trivial lines
\item the pull-backs via $\phi_n$ of the conics $\delta_2c(c-\delta_2)x^2_0 - (c-\delta_2)x^2_1 + cx^2_2=0$ for $c\ne 0,\delta_2$.
\end{enumerate}
Moreover, these curves are nonsingular and the only of them with genus $\le 2^{n-3}$ are the trivial lines, with genus $0$.
\end{lemma}
\begin{proof} Let $Y\subset X_n$ be a $\omega'_n-$integral curve. Write $Z=\phi_n(Y)$ and $Y'=\phi^*_n(Z)$. Note that $Z$ is $\omega_2-$integral, hence we have 3 cases by Lemma \ref{2integral}.\\
Suppose $Z=\{x_j=0\}\subset X_2$ is a coordinate axe. $Y'=X_n\cap\{x_j=0\}$ is nonsingular by a verification similar to the one done for $C_n$. Since that $Z$ meets all the curves $\phi(C_i)$ for $i=3,\ldots,n$ and they for the branch divisor, we have that $Y'$ is connected. Hence $Y'=Y$ and $Y$ is nonsingular. Note that $\phi_n|_Y:Y\rightarrow Z$ has degree $2^{n-2}$ and is ramified at $2^{n-2}(n-2)$ points, hence $g(Y)=2^{n-3}(n-4)+1$ by the Hurwitz formula.\\
Now suppose $Z$ is a trivial line in $X_2$. Replacing the value of $x_2$ in terms of $x_1$ in the defining equations of $X_n$ we obtain that $Y$ is a trivial line, with genus $0$.\\
Finally suppose $Z$ is a curve of type (3) in Lemma \ref{2integral}. By the same argument as in the first case, $Y'$ is connected. One can show that $Y'$ is nonsingular by a direct computation (on the affine chart $U_0$ we add the row $((c-\delta_2)x_1,-cx_2,0,\ldots,0)$ in \ref{matrix}, and for $U_1,U_2$ the computation is similar) therefore $Y=Y'$. Consider the map $\psi_n=\phi_n|_Y:Y\rightarrow Z$. If $Y$ lies above one of the curves $C_i$ then $\deg(\psi_n)=2^{n-3}$ and if $Y$ does not lie above any $C_i$ then $\deg(\psi_n)=2^{n-2}$. Anyway, $\phi_n$ is ramified at least in $(n-3)\cdot 4\cdot 2^{n-4}=2^{n-2}(n-3)$ points and $g(Z)=0$, thus for $n\ge 6$ by the Hurwitz formula we have
\begin{eqnarray*}
g(Y)>-2^{n-2}+2^{n-3}(n-3)=2^{n-3}(n-5)\ge 2^{n-3}.
\end{eqnarray*}
\end{proof}

\begin{lemma}\label{vojta} Let $\pi:X'\rightarrow X$ be a finite morphism of smooth projective surfaces over $\C$, ramified along a curve $Y\subset X'$. Let $\Lcal$ be a invertible sheaf on $X$, and take a section $\omega\in H^0(X,\Lcal\otimes S^2(\Omega^1_X))$. If $\pi(Y)$ is $\omega-$integral, then $\pi^*\omega\in H^0(X',\pi^*\Lcal\otimes S^2(\Omega^1_{X'}))$ vanishes identically on $Y$. 
\end{lemma}
\begin{proof} This is a particular case of \cite{Vojta2} Lemma 2.10.
\end{proof}

We recall to the reader that $\omega'_n=\phi_n^*\omega_2$.

\begin{lemma}\label{omegan} Define $\omega'_2=\omega_2$. The sections $\omega'_n$ determine sections 
$$
\omega_n\in H^0(X_n,\Ocal_{X_{n}}(7-n)\otimes S^2(\Omega^1_{X_n}))
$$
such that each $\omega_n-$integral curve is a $\omega'_n-$integral curve. Moreover, the $\omega_n-$integral curves are the same as the $\omega'_n-$integral curves, with the only possible exception of $\omega'-$integral curves lying over $C_3,\ldots,C_n$. 
\end{lemma}
\begin{proof} By induction. The case $n=2$ is clear. Assume it for $n=m-1$ with $m>2$. Note that $\pi_m(C_m)$ does not lie over any of the curves $C_3,\ldots,C_{m-1}$ because they have different images in $X_2$, hence $\pi_m(C_m)$ is $\omega_{m-1}-$integral by Lemma \ref{nintegral} and induction hypothesis. Consider the section 
$$
\pi^*_m\omega_{m-1}\in H^0(X_m,\pi_m^*\Ocal_{X_{m-1}}(7-(m-1))\otimes S^2(\Omega^1_{X_{m}}))=H^0(X_m,\Ocal_{X_{m}}(7-(m-1))\otimes S^2(\Omega^1_{X_{m}}))
$$ 
(recall that $\pi_n^*\Ocal_{X_{n-1}}(1)=\Ocal_{X_n}(1)$). By Lemma \ref{vojta} we have that $\pi^*_m\omega_{m-1}$ vanishes identically on $C_m$, thus $\pi^*_m\omega_{m-1}$ determines a global section $\omega_m$ in $\Ocal_{X_{m}}(7-m)\otimes S^2(\Omega^1_{X_{m}})$ by taking 
$$
\omega_m=\frac{1}{x_m}\pi^*_m\omega_{m-1}.
$$ 
Call $U_m$ the open set of $X_m$ obtained by deleting the curves lying over any of the $C_3,\ldots,C_m$. The sections $\omega'_m$ and $\omega_m$ agree on $U_m$ up to a non-vanishing factor, therefore the $\omega'_m-$integral curves and the $\omega_m-$integral curves are the same on $U_m$. A curve lying over some $C_i$ is of type (3) in Lemma \ref{nintegral} (see Equation \ref{CninP2}), hence it is $\omega'_m$-integral, and we are done.
\end{proof}
\begin{corollary}\label{lastcurves} For $n\ge 6$, the only $\omega_n-$integral curves with genus $\le 2^{n-3}$ on $X_n$ are the trivial lines, with genus $0$. 
\end{corollary}
\begin{proof} This follows from Lemma \ref{nintegral} and Lemma \ref{omegan}.
\end{proof}
\begin{theorem}\label{spset} For $n\ge 8$, the only curves of genus $0$ or $1$ on $X_n$ are the trivial lines.
\end{theorem}
\begin{proof} Let $Y\subset X_n$ be a curve of genus $0$ or $1$ and write $i:\tilde{Y}\rightarrow Y$ for its normalization. On the one hand, the curve $\tilde{Y}$ has genus $0$ or $1$, hence $\Kcal_{\tilde{Y}}$ has non-positive degree. On the other hand, the sheaf $i^*\Ocal_{X_{n}}(7-n)$ has negative degree because $n\ge 8$. Therefore, $i^*\Ocal_{X_{n}}(7-n)\otimes\Kcal_{\tilde{Y}}^{\otimes 2}$ has no nonzero global section on $\tilde{Y}$, hence $i^*\omega_n$ vanishes identically on $\tilde{Y}$. From this we deduce that $Y$ is a $\omega_n-$integral curve with genus $\le 1$ on $X_n$, and we are done by Corollary \ref{lastcurves}.
\end{proof}


\section{Proofs of results related to number fields}\label{SecNum}

We understand that, given a sequence $\delta_2,\delta_3,\ldots$ of distinct non-zero elements in $K/\Q$, the surfaces $X_n$ are defined by Equation \eqref{defining}.  

\begin{lemma}\label{correspondence} Fix a sequence $(a_1,a_2,\ldots a_n)$ in $K/\Q$, with $n\ge 3$ and pairwise distinct $a_i$. Set $\delta_i=a_i-a_1$ for $i\ge 2$. There is a bijective correspondence between the set of monic polynomials $f\in K[x]$ of degree two satisfying that $f(a_i)$ is a square for $i=1,\ldots, n$, and $X_n(K)\cap\{x_0\ne 0\}$. This correspondence is given by the map $j(f)=[1:\sqrt{f(a_1)}:\cdots :\sqrt{f(a_n)}]$ and has the property that $f$ is a square in $K[x]$ if and only if $j(f)$ lies in a trivial line of $X_n$.
\end{lemma}
\begin{proof}
Take a polynomial $f=x^2+ax+b\in K[x]$ with the property that $f(a_1)=b_1^2,\ldots,f(a_n)=b_n^2$ are squares in $K$, then 
\begin{eqnarray*}\label{verificacion}
\delta_2 b_i^2&=&(a_2-a_1)f(a_i)=(a_2-a_1)(a_i^2+ua_i+v)\\
&=&(a_i-a_1)(a_2-a_1)(a_i-a_2)\cdot 1 - (a_i-a_2)(a_1^2+ua_1+v)+(a_i-a_1)(a_2^2+ua_2+v)\\
&=&\delta_i\delta_2(\delta_i-\delta_2)1^2-(\delta_i-\delta_2)b_1^2+\delta_ib_2^2
\end{eqnarray*}
Therefore, for each polynomial $f=x^2+ux+v\in K[x]$ with the property that $f(a_1),\ldots,f(a_n)$ are squares in $K$, we have that $j(f)\in X_n(K)\cap\{x_0\ne 0\}$.\\
Conversely, given a point $p=[1:b_1:\cdots:b_n]\in X_n(K)\cap\{x_0\ne 0\}$, define 
$$
f_p=x^2+ \frac{b_2^2-b_1^2-a_2^2+a_1^2}{a_2-a_1}x+\frac{a_1a_2(a_2-a_1)-a_1b_2^2+a_2b_1^2}{a_2-a_1}\in K[x] 
$$
The polynomial $f_p$ is the only monic polynomial of degree two satisfying $f_p(a_1)=b_1^2$ and $f_p(a_2)=b_2^2$. Moreover, after a standard computation we get
$$
\delta_2 f_p(a_1+\delta_i)= \delta_i\delta_2(\delta_i-\delta_2)-(\delta_i-\delta_2)b_1^2+\delta_ib_2^2
$$
and, since $p\in X_n(K)\cap\{x_0\ne 0\}$, we obtain $\delta_2f_p(a_1+\delta_i)=\delta_2 b_i^2$. Therefore $f_p(a_i)=b_i^2$ for each $i$. \\
Clearly $j$ and $p\mapsto f_p$ are inverses, hence $j$ is bijective.

Assume that $j(f)=[1:b_1:\cdots:b_n]$ lies in a trivial line for some $f=x^2+ux+v\in K[x]$. Thus we have an equation of the kind $\pm b_2-\delta_2=\pm b_1$, say $\epsilon' b_2=\epsilon b_1 +a_2-a_1$. Therefore $b_2^2=b_1^2+2\epsilon(a_2-a_1)b_1+(a_2-a_1)^2$ and we get
\begin{eqnarray*}
\left(\frac{b_2^2-b_1^2-a_2^2+a_1^2}{a_2-a_1}\right)^2-4\left(\frac{a_1a_2(a_2-a_1)-a_1b_2^2+a_2b_1^2}{a_2-a_1}\right)=4b_1^2(\epsilon^2-1)=0
\end{eqnarray*}
So we have $f=f_{j(f)}=\left(x+\frac{u}{2}\right)^2$.
\end{proof}\\

First we prove Theorem \ref{RepNum}.\\

\begin{proof} We follow the notation of Section \ref{geometric}. For $i=2,\ldots,8$ set $\delta_i=a_i-a_1$ and note that $X_2,\ldots,X_8$ are defined over $K$. If Conjecture \ref{ConjBombieri} holds then there exists a proper Zariski closed subset $Z\subset X_8$ such that all the $K-$rational points of $X_8$ belong to $Z$. Given an irreducible curve $Y\subset X_n$, if $Y(K)$ is dense in $Y(\C)$ then $Y$ is defined over $K$ and, by Faltings' Theorem, $Y$ has genus at most $1$. Therefore we can take $Z$ as the union of a finite number of curves on $X_8$ with genus $0$ or $1$, up to a finite number of $K$-rational points. 

By Theorem \ref{spset} and Lemma \ref{correspondence} we can conclude. 
\end{proof}\\

Now we prove Corollary \ref{buchistrong}.\\

\begin{proof} Since the set $E(\Q,(a_i)_i)$ is finite, it is enough to show that a monic polynomial $f\in\Z[z]$ which is not a square, satisfies that $f(n)$ is a square at most for a finite number of $n\in \Z$. Indeed, the graph of $y=\sqrt{f(x)}$ is asymptotic to the graph of $y=|x|$ hence for large enough $|x|$ it has no integer point.
\end{proof}\\

Finally, here we have the proof of Theorem \ref{LogicNum}.\\

\begin{proof} Let $X$ be a B\"uchi surface of length $n\ge 8$ such that Bombieri's conjecture holds for $X$. Complete the sequence $\delta_2,\ldots,\delta_n$ to an infinite sequence $(\delta_i)_{i\ge 1}$ of non-zero distinct integers, and consider the corresponding surfaces $X_i$, where $X=X_n$. By an obvious modification of Corollary \ref{buchistrong}, Lemma \ref{correspondence} allows us to show that there exists an integer $M\ge n$ such that any point $[1:b_1:\cdots:b_M]\in X_M$ with $b_i\in\Z$ must lie in a trivial line. Therefore we can write a $\Lcal_2$-formula $\psi$ with the property that, $\Z\vDash \psi(c_1,\ldots,c_n)$ if and only if the $c_i$ are integer squares and $p=[1:\sqrt{c_1}:\cdots:\sqrt{c_n}]\in X_M$ lies in a trivial line. By Lemma \ref{correspondence} there exists $\nu$ such that $c_2=\nu^2$ and $c_i=(\nu+\delta_i)^2$ for $i\ge 3$. This proves the non-trivial implication of the fact that the $\Lcal_2$-formula $\Psi(x,y)$
$$
\exists c_1,\ldots,c_n (\psi(c_1,\ldots,c_n)\wedge c_2-c_1=2\delta_2 x+\delta_2^2\wedge y=c_1
$$
defines the relation $y=x^2$ in $\Z$. From here it is clear that multiplication is positive existentially definable in $\Z$ over $\Lcal_2$.
\end{proof}

\section{Proof of Theorem \ref{RepFun} (Function Fields)}\label{SecFun}

We use the same notation as in Section \ref{geometric}.

\begin{proposition} Let $n\ge 8$. If $Y\subseteq X_n$ is a curve, its normalization is $i:\tilde{Y}\to Y$ and $g(\tilde{Y})< \frac{n-3}{4}$, then $Y$ is a $\omega_n$-integral curve.
\end{proposition}
\begin{proof} Let $i:\tilde{Y}\to Y$ be the normalization map. We have
$$
i^*\omega_n\in H^0(X_n,i^*\Ocal(7-n)\otimes\Kcal_{\tilde{Y}}^{\otimes 2}).
$$
As $\deg i^*\Ocal_{X_n}(1)\ge 1$, for $n\ge 8$ we get
\begin{eqnarray*}
\deg \left(i^*\Ocal_{X_n}(7-n)\otimes\Kcal_{\tilde{Y}}^{\otimes 2}\right)=(7-n)\deg i^*\Ocal_{X_n}(1) + 4g(\tilde{y})-4\\
\le 7-n+4g(\tilde{Y})-4=4g(\tilde{Y}) +3 -n<0.
\end{eqnarray*}
therefore $i^*\omega_n$ is zero in $\tilde{Y}$.
\end{proof}\\

Now we present the proof of Theorem \ref{RepFun}.\\

\begin{proof} We can assume $F=\C$. Suppose $P$ has some non-constant coefficient and $P(a_i)=h_i^2, i=1,\ldots,M$ for some $a_i\in \C$ and $h_i\in K(C)$. Using Lemma \ref{correspondence} twice, with $K=K(C)$ and $K=\C$, one can check that $h=[1:h_1:\ldots:h_M]$ defines a non-constant morphism $h:C\to X_M$, where we consider $\delta_i=a_i-a_1$ in the definition of $X_M$. Since $C$ is a complete variety we obtain that $\im(h)$ is algebraic. Let $Y$ be an irreducible curve containing $\im(h)$, since $h$ is dominant on $Y$, we conclude that $h$ factors through $\tilde{Y}$. By Riemann-Hurwitz Formula $g(\tilde{Y})\le g(C)\le \frac{M}{4}-1<\frac{M-3}{4}$ hence $Y$ is a $\omega_M$ integral curve by the previous lemma. Therefore $Y$ is nonsingular and $g(Y)\le \frac{M}{4}-1< 2^{M-3}$ with $M\ge 8$, thus $Y$ is a trivial line by Lemma \ref{spset}. This implies that $\im(h)$ is contained in a trivial line, and the conclusion follows from Lemma \ref{correspondence}.
\end{proof}


\section{Proof of Theorem \ref{LogicMer}}\label{SecLogicMer}

We will use the positive answer of $\BP(\Mcal_p)$ and $\BP(\Acal_p)$.

Let $R$ be the ring $\Acal_p$ or the field $\Mcal_p$. The following formula  
$$	
F[x,y]: \exists u_1\cdots\exists u_{35}\left(\wedge_{i=1}^{35}P_2(u_i)\right)\wedge \left(\wedge_{i=2}^{34} u_{i-1}+u_{i+1}=2u_i+2\right) \wedge x=u_1 \wedge 2y+1=u_2-u_1
$$
is satisfied in $R$ if and only if $y=x^2$ or $x,y\in\C_p$. Then the $\Lcal_2'$-formula (actually we should use $f_z$)
$$
G[x,y]: F[x,y]\wedge F[zx,z^2y]
$$
is satisfied in $R$ if and only if $y=x^2$. Therefore, the $\Lcal_2'$-formula
$$
H[x,y,w]: \exists u\exists v \left(G[x+y,u]\wedge G[x-y,v]\wedge u=v+4w\right)
$$
is satisfied in $R$ if and only if $w=xy$. This proves Theorem \ref{LogicMer}







\end{document}